\documentclass[a4paper,12pt,twoside]{amsart}
\usepackage{amsmath}
\usepackage{amsthm}
\usepackage{amssymb}
\usepackage{fancyheadings}
\usepackage{exscale}
\usepackage{latexsym}

\newtheorem{theorem}[equation]{Theorem}
 
 \newtheorem{prop}[equation]{Proposition}
 
 \newtheorem{corr}[equation]{Corollary}
\newtheorem{cor}[equation]{Corollary}
 \newtheorem{lemma}[equation]{Lemma}
 \theoremstyle{definition}

 \newtheorem{remark}[equation]{Remark}
 
 \newtheorem{example}[equation]{Example}

 \DeclareMathOperator{\Gl}{Gl}
 \DeclareMathOperator{\Sl}{Sl}

 \DeclareMathOperator{\charo}{char}

\DeclareMathOperator{\Quot}{Quot}

 \newcommand{\field}[1] {\mathbb #1}

 \newcommand{\F} {{\field F\:\!}}
 \newcommand{\K} {\field K\:\!}

 \newcommand{\A}{{\mathcal A}}
\newcommand{\sR}{{\mathcal R}}

\newcommand{\om}{\omega} 

\newcommand{\highexpo}{\,\, \rule[0ex]{0ex}{0ex}m_{\hspace{-.2ex}H\rule[1ex]{0ex}{0ex}} }
\newcommand{\expo}{ \rule[0ex]{0ex}{2ex}m_{\hspace{-.1ex}\rule[0ex]{0ex}{1.5ex}H}}

\newcommand{\del}{\partial}

 \def \verylongrightarrow{\hbox to 28pt{\rightarrowfill}}
\def \verylongleftarrow{\hbox to 28pt{\leftarrowfill}}
 \def \elra{\hbox to 30pt{\rightarrowfill}}
\def \ella{\hbox to 30pt{\leftarrowfill}}

\begin{document}

\title[Jacobians of reflection groups]
{Jacobians of reflection groups}

\author[Julia Hartmann]{Julia Hartmann}
\email{Julia.Hartmann@iwr.uni-heidelberg.de}
\address{IWR, Universit\"at Heidelberg, 69120 Heidelberg, Germany}

\author[Anne V.\ Shepler]{Anne V.\ Shepler}
\email{ashepler@unt.edu}
\address{Dept.\ of Mathematics,
         University of North Texas,
         P.O. Box 311430,
         Denton, TX, 76203}


\thanks{Work of second author partially supported by National Security Agency 
grant MDA904-03-1-0005}

\begin{abstract}
Steinberg showed that when a finite reflection group acts on a real or complex vector space of finite dimension,
the Jacobian determinant of a set of
basic invariants 
factors into linear forms which define the reflecting hyperplanes.
This result generalizes verbatim to fields whose characteristic is 
prime to the order of the group. 
Our main theorem gives a generalization of Steinberg's result for arbitrary
fields using a ramification formula of Benson and Crawley-Boevey.
As an intermediate result, we show that every finite group which 
fixes a hyperplane pointwise has a polynomial ring of invariants.
\end{abstract}

\hphantom{x}\vspace{-3ex}
\maketitle

\section{Introduction}
\label{introduction}
The advent of modern algebra owes much to invariant theory.
Many of our classical theorems arose from Noether's investigations of a finite group
$G \leq \Gl(V)$ acting linearly on an $n$-dimensional vector space
$V$ over a field $\F$.  
The action of $G$ induces a natural action on the
 polynomial ring $\F[V]\cong \text{Sym}(V^*)$. 
Noether showed that the ring $\F[V]^G$ of invariant polynomials is 
finitely generated as an algebra. 
We are interested in the case when generators of $\F[V]^G$ are 
algebraically independent:  
we say that $G$ has a {\bf polynomial ring of invariants}
if $\F[V]^G = \F[f_1, \ldots, f_n]$ for some homogeneous
polynomials $f_i$ called {\bf basic invariants}
or 
{\bf fundamental invariants}.
Although there are many choices of basic invariants, 
their degrees are unique, and thus the integers 
$\deg f_1 - 1, \ldots, \deg f_n -1$
depend only on the group.
We call these integers the {\bf exponents} of $G$.
When $G$ has a polynomial ring of invariants, we 
define the {\bf Jacobian determinant} 
$J=J(f_1, \ldots, f_n)=  \det (\del f_i/\del z_j)$. This polynomial is 
nonzero (see \cite[5.4]{benson}) and well-defined up to a
nonzero element of $\F$ depending on the choice of basic invariants and 
basis  $\{z_j\}$ of $V^*$.
In this article, we examine the structure of the Jacobian determinant $J$.
No assumption is made on the ground field $\F$.

Elements of finite order in $\Gl(V)$ which
fix a hyperplane pointwise are called {\bf reflections}.
For any subgroup $G$ of $\Gl(V)$ and hyperplane $H$ in $V$, 
define 
$$G_H=\{g \in G:\; g{\vert}_{{\rule[0ex]{0ex}{1.5ex}}_H}=\operatorname{id}_H\},$$ 
the pointwise stabilizer of $H$ in $G$. The hyperplanes $H$ for which
$G_H$ is nontrivial are called {\bf reflecting hyperplanes} of $G$.  
For each hyperplane $H$ in $V$, let $l_H$ in $V^*$
be a linear form with $\ker l_H = H$.

In Section~2, we examine 
the set of root vectors of reflections about a common hyperplane.
We show that every finite subgroup of $\Gl(V)$ which fixes
a hyperplane $H$ in $V$ pointwise has
a polynomial ring of invariants, and furthermore,
that the Jacobian determinant is a power of $l_H$. 

Our main result (Theorem~\ref{main}) states that if $G\leq\Gl(V)$ is a finite
group with a polynomial ring of invariants, its Jacobian determinant
factors as a product of
linear forms defining the reflecting hyperplanes of $G$. The multiplicity with
which each linear form $l_H$ occurs is the sum of the exponents of $G_H$.


This theorem has roots in the rich theory of reflection groups.
A finite subgroup of  $\Gl(V)$ is a 
{\bf reflection group}
if it is generated by reflections,
and the collection $\A$ of its reflecting hyperplanes
is called the {\bf reflection arrangement}.
When the characteristic of $\F$ is prime to $|G|$, 
a well-known theorem of Serre and Shephard, Todd, and Chevalley
(see \cite[Ch.~7]{Poly}) 
states that a finite subgroup of $\Gl(V)$ is a reflection group if and only if 
it has a polynomial ring of invariants.
Steinberg~\cite{steinberg60} 
showed in this case that the Jacobian determinant of a set of basic
invariants factors into powers of linear forms defining the reflecting hyperplanes:
$$
J  \ \doteq \ \prod_{H \in \A} l_H^{\,|G_H|-1}
$$
(we write $a \doteq b$
to indicate that $a$ and $b$ are equal up to a nonzero constant).
In particular, 
the above factorization holds for all Weyl groups, Coxeter groups,
and complex reflection groups (see \cite[Thm.~6.42]{orlikterao}).
In this case, each $G_H$ is a cyclic group and the
only nonzero exponent of $G_H$ is 
$|G_H|-1$.

Serre~\cite{serre} showed that in arbitrary characteristic, 
every finite subgroup of $\Gl(V)$ 
with a polynomial ring of invariants must be generated by reflections. 
The converse may fail when the characteristic of $\F$ divides the
order of $G$ 
(for example, see \cite{kemper97}).
Unfortunately, Steinberg's characterization 
of the Jacobian determinant in terms of the integers $|G_H|$ no longer holds
over arbitrary fields. 
The stabilizer subgroups 
$G_H$ may not be cyclic, in which case the integers $|G_H|-1$ as $H$ 
runs over all reflecting hyperplanes 
will usually not sum to the degree of $J$.
This is no surprise, 
as the class of reflections
is larger in some sense over an arbitrary field
than over 
a characteristic zero field. 
The reflections in $\Gl(V)$ not only include diagonalizable
reflections (with a single nonidentity eigenvalue), but also
{\bf transvections}, reflections with determinant~$1$ which can
not be diagonalized. 
The transvections in $\Gl(V)$ prevent one from developing a theory
of reflection groups mirroring that for Coxeter groups or 
complex reflection groups. 
(For example, even if a reflection group has a polynomial ring of invariants, 
the Jacobian $J$ may be invariant or lie in the Hilbert ideal
generated by the basic invariants---see 
Section~\ref{examples}.)
If $G$ lacks transvections, then it shares some characteristics
with reflection groups over characteristic zero fields, for example,
the pointwise stabilizer of any hyperplane in $V$ is cyclic. One can
deduce that Steinberg's description remains valid in this special case (see
\cite{hartmann}).

Theorem~\ref{main} implies a conjecture by Victor Reiner:
if $G$ has a polynomial ring of invariants,
then the zero locus of the Jacobian determinant is exactly the union of the reflecting hyperplanes.
Reiner, Stanton, and Webb~\cite{reiner04}
use this corollary in generalizing
Springer's theory of regular numbers
in characteristic zero to arbitrary fields.

\section{ONE HYPERPLANE}

In this section, we consider finite groups $G\leq \Gl(V)$ 
that fix a hyperplane $H$ in $V$ pointwise. Do such groups have polynomial
invariants?
Smith~\cite[Chap.~8]{Poly} examines these groups following Landweber and Stong
(\cite{land}), who described a set of basic invariants over $\F_p$ for $p$ a prime.  
Unfortunately, parts of the description 
(end of Section 2 in \cite{Poly}) do not extend to other finite fields
(see Example~\ref{ex1} in the last section).
Nakajima~\cite{nakajima} proves that if the field
$\F$ is finite and if every nonidentity element of $G$ is a transvection,
then $G$ has a polynomial ring of invariants.
Benson~\cite{benson} (after Thm~3.13.2) 
indicates that the situation may be quite different over arbitrary fields.
(Note that 
$G$ may not be conjugate to a group with entries in a 
finite field.)
We show below that any finite subgroup of $\Gl(V)$
which fixes one hyperplane pointwise 
must indeed have a polynomial ring of invariants (no assumptions on the field).
We begin by examining the structure of such groups.  

Let $H$ be a hyperplane in $V$ defined by some linear form
$l_H \in V^*$.
For any reflection $t \in\Gl(V)$ which fixes $H$ 
pointwise, let $\alpha_t$ be the 
{\bf root vector}
of $t$ (with respect to $l_H$) defined by 
$$ 
 t(v) = v + l_H(v)\, \alpha_t \ \ \quad\text{for all } v \in V.
$$
Note that a transvection is a reflection whose root
vector lies in the reflecting hyperplane, i.e., $l_H(\alpha_t) =0$
(see for example
\cite{smithneusel}, section~6.2).
For any set $S$ of reflections, let $\sR(S)$ be the corresponding
set of root vectors in $V$.

If char $\F=0$, then any group $G$ which fixes a 
hyperplane pointwise is necessarily cyclic 
and its order equals the maximal order of
a diagonalizable reflection in $G$ (which then generates $G$). 
The next lemma gives the order of $G$ when char $\F >0$.

\begin{lemma}
\label{order}
Assume that $\charo(\F)=p>0$.
Suppose $G \leq \Gl(V)$ is a finite group which fixes a hyperplane $H$ in $V$
pointwise.
Let $\sigma$ be a diagonalizable reflection in $G$ 
of maximal order with 
$\om = \det(\sigma)$.
Let $K\trianglelefteq G$ be the 
subgroup generated by the transvections in $G$. Then
\begin{enumerate}
\item 
The action of $\sigma$ on $K$ by conjugation
translates into multiplication by
$\omega$ on $\sR(K)$ and thereby endows
$\sR(K)$ with the structure of an $\F_p(\om)$-vector space.
\item
$T\subset K$ is a minimal set satisfying $G=\langle T,\sigma\rangle$
if and only if $\sR(T) \subset \sR(K)$
is a basis for $\sR(K)$ over $\F_p(\om)$.
\item 
The group $G$ has order
$|\om| \cdot |\F_p(\om)|^d$, where $d$ is the minimal
number of transvections needed to generate $G$ together with $\sigma$.
\end{enumerate}
\end{lemma}
\begin{proof}
If $s$ and $r$ are both reflections about $H$, 
then the root vector of the product is a
linear combination of the root vectors:
\begin{equation}
\label{star}\tag{$\ast$}
 \alpha_{sr}\ =\ \lambda_{r}\, \alpha_{s} + \alpha_r,
\end{equation}
where $\lambda_r = 1 + l_H(\alpha_r)$ 
is the nonidentity eigenvalue of $r$
(if $r$ is not a transvection)
or the eigenvalue $1$ (if $r$ is a transvection).
In particular,
$\alpha_{st} = \alpha_s + \alpha_t$ for all $s\in G$, $t \in K$.
Note that $\alpha_{\sigma} = -\om \alpha_{\sigma^{-1}}$.  
Fix some transvection $t$ in $K$.
An easy computation then shows that
 \begin{equation}\label{starstar}\tag{$\ast\ast$}
   \alpha_{\sigma^{-1} t \sigma}\ =   \ \om \alpha_t \quad\text{and
    thus}\quad 
    \alpha_{\sigma^{-1} t^m \sigma}\ = \ \om \alpha_{t^m}\ =  \   m\om\alpha_t
 \end{equation}
for $m \in \{0,1,\ldots,p-1\}$.

We claim that the root vector of any element 
in the subgroup $\langle\sigma, t \rangle$ 
must lie in the 
$\F_p(\omega)$-span 
of $\alpha_\sigma$ and $\alpha_t$. Indeed, 
we can write the element as a product of the generators $\sigma$ and $t$
and use Equation~(\ref{star}) repeatedly.
In particular, the root vector of a transvection 
in $\langle\sigma, t \rangle$ 
lies on $H$ and thus must be an $\F_p(\omega)$-multiple of $\alpha_t$ alone 
(as $\alpha_t$ lies on $H$ but $\alpha_\sigma$ does not).
On the other hand, any $\F_p(\omega)$-multiple of $\alpha_t$
is the root vector of some transvection in $\langle\sigma, t\rangle$:
if 
$e=|\om|$ and 
$$
 \alpha = (m_1\om+m_2\om^2+\ldots+m_e\om^e)\
       \alpha_t \quad \text{ for some } m_j \in \F_p,
$$
then $\alpha$ is the root vector of the transvection
$$
  (\sigma^{1}t^{m_1}\sigma^{-1})
       (\sigma^{2}t^{m_2}\sigma^{-2})
       \cdots
       (\sigma^{e}t^{m_e}\sigma^{-e})
$$
in $\langle\sigma, t\rangle$ by 
Equation~(\ref{starstar}).
Since each transvection about $H$ is determined by its root vector, 
the transvections in 
$\langle\sigma,t\rangle$ correspond bijectively
to the $\F_p(\om)$-multiples of $\alpha_t$.

More generally, if $t_1, \ldots, t_k$ are transvections in $G$, 
then a similar argument shows that the
$\F_p(\om)$-span of $\alpha_{t_1},\ldots,\alpha_{t_k}$ 
is the set of root vectors corresponding to the group
$K  \cap  \langle\sigma, t_{1},\ldots t_{k}\rangle$.
This proves part~(1).

Suppose $T=\{t_1, \ldots, t_k\}$ is a minimal subset of $K$ satisfying
$G=\langle T,\sigma\rangle$ and let $\alpha_i=\alpha_{t_i}$.
Then no $t_m$ lies in the group generated by $\sigma$ and  $\{t_i:i\neq m\}$
and hence no $\alpha_m$ is an $\F_p(\om)$-linear combination
of $\{\alpha_i:i\neq m\}$. 
Thus, the root vectors $\alpha_1, \ldots, \alpha_k$
are linearly independent over $\F_p(\om)$.
As $T$ generates $G$ together with $\sigma$, 
the root vectors $\alpha_1, \ldots, \alpha_k$ span $\sR(K)$ over $\F_p(\om)$
and hence form a basis. 
Conversely, if the root vectors of some set $T$
form an $\F_p(\om)$-basis of $\sR(K)$, then $T$ is a minimal
subset of $K$ generating $G$ together with $\sigma$, which proves~(2).

Finally, if $\sR(T)$ is a basis of $\sR(K)$ over $\F_p(\om)$
for some $T \subset K$, then 
$|K| = |\sR(K)|=
|\operatorname{span}_{\F_p(\om)}\sR(T)| 
= |\F_p(\om)|^{|T|}$,
and hence $|G| = |\sigma| \cdot |K|=|\om|\cdot |\F_p(\om)|^{|T|}$.
\end{proof}

Recall that a polynomial 
$f \in \K[x_1,\ldots,x_r]$ 
over a field $\K$ is called {\bf
  additive} if it induces an additive homomorphism $\K^r \rightarrow \K$. In
characteristic zero, the only additive polynomials are linear polynomials without constant term.
If the characteristic of $\K$ is $p>0$ and $|\K|=\infty$, then additive polynomials
are exactly those in which all exponents are $p$\,-powers
(for example, see \cite{lang}, VI~\S 12). 
(In particular, every additive
polynomial is ${\mathbb P}$-linear for the prime field ${\mathbb P}$ of $\K$.)

\begin{lemma}\label{additive}
Let $\F$ be a field and let $A\subseteq \F$ be a finite additive subgroup.
Then the polynomial
$f(X)=\prod\limits_{a \in A} (X+a)\in\F[X]$ is additive.
\end{lemma}
\begin{proof}
Consider the polynomial 
$$
F(X,t) := f(X+t) - f(X) = a_{m-1}(t) X^{m-1} + \ldots + a_1(t) X + a_0 (t),
$$ where $t$ is
another variable and the $a_i$ are polynomials in $t$. Note that $a_0(t) =
F(0,t) = f(t)$ by definition, and $\operatorname{deg}_t(a_i)<m$ for $i\geq 1$.  
But $F(X,t_0)=0$ for all $t_0 \in
A$, since $f(X+t_0) = f(X)$ by definition of $f$.
Hence, for every $t_0 \in A$, the coefficients 
$a_1(t_0),\ldots,a_{m-1}(t_0)$ are all zero. Thus for $i\geq 1$,
the polynomial $a_i(t)$ has at least $m$ zeroes. 
Since each $a_i$ has degree at most
$m-1$ in $t$, it must be identically zero. 
This shows that
$f(X+t)-f(X) = F(X,t)=a_0(t) = f(t)$, and so $f$ is additive. 
\end{proof}
The reader who is familiar with the Landweber-Stong invariants over prime fields
is encouraged to peruse Example~\ref{ex1} in the last section before
considering the next proposition and its proof.
\begin{prop}
\label{hyper}
Let $H\leq V$ be a hyperplane defined by some $l_H\in V^*$.
Let $G \leq \Gl(V)$ be a finite group
fixing $H$ pointwise. Then
\begin{enumerate}
\item\label{a} The group $G$ has a polynomial ring of invariants.
\item\label{b} 
The Jacobian determinant is 
$
 J\ \doteq \ l_H^{\, m\rule[-.3ex]{0ex}{0ex}},
$
where $m$ is the sum of the exponents of $G$.
\end{enumerate}
\end{prop}
\begin{proof}
The group $G$ is generated by a diagonalizable reflection $\sigma$
with eigenvalue $\omega$ of order $e$ together
with a minimal set of transvections $t_1, \ldots, t_r$
(see Lemma~\ref{order}).
For $k=1,\ldots,r$, let $G_k = \left<\sigma, t_1, \ldots, t_k\right>$, 
and let $G_0=\left<\sigma\right>$.
Choose a basis $e_1,\ldots,e_n$ of $V$ such that $\sigma$ is in diagonal form
and $e_1,\ldots,e_{n-1}$ span $H$. Let $z_1,\ldots,z_n$ be the dual basis of
$V^*$ and rescale $l_H$ so that $z_n=l_H$.
Consider the case $\charo(\F)=p>0$.

We prove by induction on $k$ a stronger statement:
$\F[V]^G=\F[f_1,\ldots,f_n]$ for some homogeneous polynomials $f_i$ where
$f_n=z_n^e$, $J(f_1, \ldots, f_n)=z_n^{m}$,
and for $i<n$, the degree of $f_i$ is a $p$\,-power and $f_i$
is additive as a polynomial in  $\F(z_n)[z_1,\ldots,z_{n-1}]$.
For $k=0$, these claims are satisfied by setting
$f_n=z_n^e$ and $f_i = z_i$ for $i<n$.
Note that these are also the basic invariants when the
characteristic of $\F$ is zero (as $G=\langle \sigma \rangle$ in this case).

Let $k\geq 0$ and 
assume the induction hypothesis holds for the group $G_k$
with $\F[V]^{G_k} = \F[f_1, \ldots,f_n]$.
Let $d_i$ be the degree of each $f_i$ and
let $t=t_{k+1}$.
By our choice of basis, 
$t(z_i) = z_i + a_i z_n$ for some $a_i\in\F$ when $i<n$ and $t(z_n)=z_n$.
For $i<n$, each $f_i$ is additive over the infinite field $\F(z_n)$, and thus 
$$
\begin{aligned}
 t f_i(z_1, \ldots, z_n) &= f_i(z_1 + a_1 z_n, \ldots,z_{n-1}+a_{n-1}z_n,z_n)\\
 & = f_i(z_1, \ldots, z_n) + f_i(a_1 z_n, \ldots, a_{n-1}z_n, z_n)
 = f_i + b_i z_n^{d_i}
\end{aligned}
$$
for some $b_i \in \F$ (note that the second summand only depends on the variable $z_n$).
Thus $b_i=0$ exactly when $f_i$ is invariant under $t$.

Relabel $f_1, \ldots, f_{n-1}$ so that $f_1$ has minimal degree
among those $f_i$ which are not invariant under $t$.
Define $f'_2, \ldots, f'_{n-1}$ by
$f_i'= f_i + c_i f_1^{d_i/d_1}$
where $c_i =  -b_i/b_1^{d_i/d_1}$.
The constants $c_i$ are chosen so that each $f'_i$ is invariant under $t$.
The degrees of $f_2', \ldots, f'_{n-1}$
are again $p$\,-powers, since
$d_i/d_1$ is a nonnegative $p$\,-power whenever
$c_i \neq 0$.
Furthermore, $f'_2,\ldots, f'_{n-1}$ 
are additive over $\F(z_n)$ as they are the compositions of 
additive homomorphisms.  
Define $f_n'=f_n$.
Then $f'_2, \ldots, f'_{n}$
are invariant 
under $t$ and $\sigma,t_1, \ldots, t_k$ (as each $f_i \in \F[V]^{G_k}$).
Hence, $f'_2, \ldots, f'_n$ are invariant under $G_{k+1}$.

We take the product over the orbit of  $f_1$ to produce a polynomial
$f'_1$ invariant under $t$.
Define 
$$
   h(X) = \prod_{a \in A} (X+ a z_n^{d_1})
               \ \ \in \F(z_n)[X],
$$
where
$A =\F_p(\omega)b_1=\{\lambda b_1|\; \lambda \in \F_p(\omega)\}$. 
By Lemma~\ref{additive},
$h(X)$ is additive as a polynomial in $\F(z_n)[X]$.
Define $f'_1 = h(f_1) \in\F[z_1, \ldots, z_n]$.
Then $f'_1$ is additive in $\F(z_n)[z_1, \ldots, z_{n-1}]$
as it is the composition of additive homomorphisms.
 The polynomial $f'_1$ is invariant under $t$ (by its very definition)
and invariant under  $t_1, \ldots, t_k$ since both $f_1$ and $z_n$ are.
The polynomial $f'_1$ is also invariant under the diagonalizable reflection
$\sigma$ 
since $\sigma(f_1) = f_1$, $\sigma(z_n) = \om^{-1} z_n$, and $A$ is closed
under multiplication by $\F_p(\om)$. 
(In particular, $f'_1$ is a polynomial in $f_1$ and $f_n$.)
Hence, $f'_1, \ldots, f'_n \in \F[V]^{G_{k+1}}$.

We consider each $f'_i$ as a polynomial in $f_1, \ldots, f_n$.
By the chain rule,
$$
 J(f'_1,\ldots, f'_n)\ =\   J(f_1, \ldots, f_n)\ 
 \det\left( \frac{\del f'_i}{\del f_j}\right).
$$
The matrix  $\left(\del f'_i/\del f_j\right)$ is upper triangular 
with determinant $\frac{\del f'_1}{\del f_1}$.
Since $h$ is additive as a polynomial
in $\F(z_n)[X]$ by Lemma~\ref{additive}, 
every exponent of $X$ in an expansion of $h$ is a $p$\,-power.
If we expand $f'_1=h(f_1)$ as polynomial in 
$f_1$ and $z_n$, every exponent of $f_1$ will thus also be a $p$\,-power.
In particular\footnote{One can use the fact that $\F_p(\om)$ is the splitting
  field of the polynomial $X^{|A|} - X$ to compute the coefficient; 
  we do not use this coefficient in what follows.},
$$
 \del f'_1/\del f_1 \ =\ -b_1^{|A|-1}\ z_n^{d_1(|A|-1)},
$$
and hence
$
 J(f'_1,\ldots,f'_n) \doteq J(f_1,\ldots,f_n)\cdot z_n^{d_1(|A|-1)}.
$
By the induction hypothesis, $J(f_1,\ldots,f_n)$ is a power of $z_n$
and the exponent of $z_n$  is the sum of the exponents of $G_k$.
Substituting this into the last equality shows that assertion~(\ref{b}) holds
for $G_{k+1}$ and thus for $G$ by induction.
The polynomials $f_1',\ldots, f_n'$ form a set of basic
invariants for $G_{k+1}$
if and only if $J(f'_1, \ldots, f'_n)$ is nonzero 
and the product of the degrees of the $f'_i$ 
is the order of the group $G_{k+1}$ 
(for example, see 
\cite[Prop.~16]{kempercalc}).
By Lemma~\ref{order} and the induction hypothesis,
$$\deg f'_1 \cdots \deg f'_n
= |\F_p(\om)| \deg f_1 \deg f_2 \cdots \deg f_n
= |\F_p(\om)| |G_k|
= |G_{k+1}|,
$$  
and ~(\ref{a}) follows.
\end{proof}
The proof of Proposition~\ref{hyper} shows an interesting fact.
The polynomials $f_1,f_2',\ldots,f_n'$ in the 
induction step of the proof form a set of basic
invariants for ${G_k}$. 
Thus, if we choose basic invariants of $G_k$ wisely,  
we need only adjust {\em one} of them to produce 
basic invariants for $G_{k+1}$:
\begin{corr}
Let $G\leq \Gl(V)$ be a finite group which fixes a hyperplane 
$H$ in $V$ pointwise.
Let $G'=\langle G,t\rangle$ where $t \notin G$ is a
transvection 
about $H$. Then there exist basic invariants $f_1,\ldots,f_n$ for $G$
and an invariant $f_1'$ for $G'$ such that $f_1',f_2,\ldots,f_n$ form a set of
basic invariants for $G'$ (in particular, all basic invariants for $G$ except
one can be chosen to be invariant under $G'$).
\end{corr}
%
%
\section{The Jacobian Factors}

In this section, we consider a finite group $G \leq \Gl(V)$ with 
a polynomial ring of invariants. 
We show that the Jacobian determinant factors into powers of linear forms
defining the reflecting hyperplanes.
We begin with an easy consequence of Proposition~\ref{hyper}:   
\begin{lemma}
\label{LL'dividesJ}
Assume that $G\leq\Gl(V)$ is a finite group with a polynomial ring of invariants.
Let $\A$ be the reflection arrangement of $G$.
Then the Jacobian determinant is divisible by
$$
 \prod_{H \in \A} l_H^{\highexpo},
$$
where each $\expo$ is the sum of the exponents of the pointwise stabilizer $G_H$.
\end{lemma}
\begin{proof} 
Fix some reflecting hyperplane $H \in \A$.
By Proposition~\ref{hyper}, $G_H$ has a polynomial
ring of invariants and 
$$
J(f^H_1, \ldots, f^H_n) \doteq  z_n^{\highexpo},
$$
where $f^H_1, \ldots, f^H_n$ are basic invariants for $G_H$. Let 
$f_1,\ldots,f_n$ denote basic invariants for $G$. 
Each $f_i$ is invariant under $G_H$ and hence may be written as a polynomial in
the $f_j^H$. Thus
$J(f^H_1, \ldots, f^H_n) = l_H^{\highexpo}$
divides
$J(f_1, \ldots, f_n)$ by the chain rule.
The claim then follows since the linear 
forms $l_H$ for different reflecting hyperplanes 
are pairwise coprime and $\F[V]$ is a unique factorization domain.
\end{proof}

We next verify that we have found all factors of $J$.
We compare degrees using the following version of the ramification formula of
Benson and Crawley-Boevey (Corollary~3.12.2 in \cite{benson}):
\begin{lemma}\label{ram}\
Assume that $G\leq\Gl(V)$ is a finite group.
Then
$$ 
 |G|\ \psi(\F[V]^G) \ =\ \sum \limits_{H \leq V} \ |G_H|\ \psi(\F[V]^{G_H})
$$
(the sum runs over all hyperplanes in $V$).
Here $\psi(M)$ denotes the coefficient of $\frac{1}{(1-t)^{n-1}}$ in the expansion at $t=1$ of the
Poincar\'e series of a finitely generated $\F[V]^G$-module $M$.
\end{lemma}

We apply the above lemma and obtain
\begin{lemma}\label{degree}
Assume that $G\leq\Gl(V)$ is a finite group with a polynomial ring of invariants
and let $\mathcal{A}$ be its reflection arrangement.
Let $J$ be the Jacobian determinant of $G$. Then
$$
 \operatorname{deg}(J)=
  \sum \limits_{H \in \mathcal{A}} \expo ,
$$
where each $\expo$ is the sum of the exponents of the pointwise stabilizer
$G_H$.
\end{lemma}
\begin{proof}
By Proposition~\ref{hyper}, each $G_H$ has a polynomial ring of 
invariants. The product of the degrees 
$d_i^{H}$ of basic invariants for $G_H$ equals the order of $G_H$ by
\cite[Prop.~16]{kempercalc}.
Consequently, 
\begin{alignat*}{2}
\frac{1}{2} \operatorname{deg}(J)\ & =\ \frac{1}{2}\sum \limits_{i=1}^n (d_i -1) 
=\ |G| \ \psi(\F[V]^G) \\
& =\ \sum \limits_H |G_H|\ \psi(\F[V]^{G_H}) &\text{\hfill by Lemma~\ref{ram}}&\\
& =\ \sum \limits_H |G_H| \ \frac{1}{2|G_H|}\ \sum \limits_{i=1}^n (d_i^{(H)}
-1)=\ \frac{1}{2}\sum \limits_H   \expo.
\end{alignat*}
which proves the claim.
\end{proof}

The main theorem is a direct consequence.
\begin{theorem}
\label{main}
Assume that $G\leq\Gl(V)$ is a finite group with a polynomial ring of invariants.
Then the Jacobian determinant $J$ factors into a
product of powers of linear forms defining the reflecting hyperplanes.  
In fact,
$$
J \ \doteq \ \ \prod_{H \in \A} \ l_H^{\highexpo},
$$
where $\expo$ denotes the sum of the exponents of the pointwise stabilizer $G_H$.
\end{theorem}
\begin{proof}
By Lemma~\ref{LL'dividesJ}, the right hand side of the equation divides $J$.
By Lemma~\ref{degree}, both sides have the same degree.
Thus they are equal up to a scalar.
\end{proof}

We immediately obtain
\begin{cor}\label{union} 
Assume that $G\leq\Gl(V)$ is a finite group with a polynomial ring of invariants.
Then the zero set of the Jacobian determinant 
is the union of all reflecting hyperplanes of $G$.
\end{cor}

\begin{remark}
There is a geometric proof of Corollary~\ref{union} suggested by
W.~Messing: Since the extension of quotient fields
$\Quot(\F[V])/\Quot(\F[V]^G)$ is separable of degree $|G|$, the associated 
quotient  morphism 
$\pi:V\rightarrow V/G$ is an \'etale covering in a neighborhood of a point $v \in
V$ if and only if $G$ acts freely on $v$. By (a variant of) a theorem of
Serre, this happens if and only if $v$ avoids all reflecting hyperplanes. 
Since $\F[V]^G$ is a polynomial ring, $\pi$ is a morphism of affine varieties,
and thus it is \'etale near $v$ if and only if the Jacobian matrix evaluated at
$v$ is invertible, i.e., has nonzero determinant. See \cite{reiner04} for the
full argument.
\end{remark}
\section{Examples}
\label{examples}
We first give an example to illustrate Lemma~\ref{order} and 
Proposition~\ref{hyper} and 
also to clarify what goes wrong with the proofs of 
Theorems 8.2.14 and 8.2.19 in \cite{Poly} when the ground field does not have prime order.
We then give examples illustrating Theorem~\ref{main}.
\setcounter{equation}{0}

\begin{example}\label{ex1}
Consider the group $G\leq\Gl_n(\F)$ generated by the matrices
$$
  A=\left(\begin{matrix}
    1 & 0 & a \\
    0 & 1 & 0\\
    0 & 0 & 1
   \end{matrix}\right),\quad
 B=\left(\begin{matrix}
    1 & 0 & 0 \\
    0 & 1 & b\\
    0 & 0 & 1
    \end{matrix}\right), \quad \text{and} \quad
 C=\left(\begin{matrix}
    1 & 0 & c \\
    0 & 1 & c \\
    0 & 0 & 1
\end{matrix}\right)
$$
where $\F=\F_p(a,b,c)$ 
for some 
nonzero $a$, $b$, $c$ in $\F$.
The group $G$ fixes the hyperplane defined by the equation $z_3=0$.
Lemma~\ref{order} is transparent in this example:
The set $\sR(G)$ of root vectors of $G$ is just the
$\F_p$-span of $(a,0,0), (0,b,0)$ and $(c,c,0)$ in $\F^3$.
The minimum number of vectors needed to span $\sR(G)$ is exactly the minimum number
of group elements needed to generate $G$.  
Thus, $G$ can be generated by $d$ elements and no fewer
if and only if the dimension of $\sR(G)$ over $\F_p$ is $d$.
 
Assume $A$, $B$, and $C$ form a minimum generating set for $G$.
(In other words, either 
$a$ and $c$ are independent over $\F_p$ or 
$b$ and $c$ are independent over $\F_p$.)
The group $G$ has a polynomial ring of invariants with
basic invariants
\begin{align*}
f_1&= (z_1^p -a^{p-1} z_1 z_3^{p-1})^p
- c^{p-1}(a^{p-1} - c^{p-1})(z_1^p -a^{p-1} z_1 z_3^{p-1}) z_3^{p(p-1)},\\
f_2 &= (z_2^p - b^{p-1} z_2 z_3^{p-1}) -
\left(\frac{b^{p-1}-c^{p-1}}{a^{p-1} - c^{p-1}}\right)
(z_1^p -a^{p-1} z_1 z_3^{p-1}),\\
f_3&=z_3.
\end{align*}
These basic invariants are given by the proof of Proposition~\ref{hyper}.

In the special case where $a=b=1$ and $c$ is algebraic over $\F_p$ but $c\notin\F_p$,
the group $G$ is defined over the finite field $\F_p(c)=\F$,
yet Theorem 8.2.14 and the proof of Theorem 8.2.19 in \cite{Poly}
do not describe the group and basic invariants.
\end{example}
\begin{example}

Let $\F$ be the finite field $\F_q$ 
and let $\A$ be the set of all hyperplanes in 
$V=\F^n$.
For each hyperplane $H$, choose some $l_H \in V^*$ with $\ker l_H = H$,
and let $Q$ be the product of all these linear forms:
$ Q = \prod_{H \in \A} l_H $.
Then $Q$ has degree 
$ 
 |\A| =  {|V|}/{|\F^*|} = ({q^{n}-1})/({q-1}).
$

The group $G = \Gl_n(\F)$ is generated by reflections about
all hyperplanes in $V$.
The invariants of $G$ form a polynomial ring:
$
\F[V]^G \ = \ \F[f_1, \ldots, f_n],
$
where $f_{i+1}$ is the Dickson polynomial
$$
d_{n,i} \ = \ \sum_{
         \substack{W \leq V\\ \rule[.7ex]{0ex}{1ex}\dim W \, = \, i}} 
         \ \ \prod_{\substack{v \, \in \, V^*,\\ v|_W\neq 0}} v
$$
of degree $q^n-q^i$
(for example, see \cite[Prop. 8.1.3]{benson}).
Note that $f_1 = d_{n,0} \ = \ Q^{q-1}$.

Fix some hyperplane $H$ in $V$ 
and let $G_H$ be its pointwise stabilizer in $G$.
In an appropriate basis (with $z_n = l_H$),
$\F[V]^{G_H} = F[u_1, ..., u_n]$
where 
$$
   u_i= z_i^q -  z_n^{q-1} z_i \ \text{ for }\  i<n\
   \text{ and }\ u_n=z_n^{q-1}
$$
(see \cite{land}).
Note that the sum of the exponents
of $G_H$ is 
$$\expo = (n-1)(q-1)+(q-2) = n(q-1) - 1.$$
Since each $f_i$ lies in $\F[V]^{G_H}$,
each $f_i = h_i(u_1, \ldots, u_n)$
for some $h_i \in \F[V]$. 
Then
$   \del f_i/\del z_k = \sum_{j} (\del h_i/\del u_j)(\del u_j /\del z_k)$
is divisible by $z_n^{q-2}$ if $k=n$
and divisible by $z_n^{q-1}$ otherwise.
Hence, $J=\det\left(\del f_i/\del z_k\right)$ is divisible by $z_n^{\highexpo}$. 
As $H$ was arbitrary, 
$\prod_{H \in \A} l_H^{\highexpo}$ divides $J$.
But one may check that $\deg J 
 =\deg Q \cdot \left(n(q-1)- 1\right)$,
and thus 
$$
 J \ =\ Q^{n(q-1)-1} \ =\ \prod_{H\in\A}l_H^{\highexpo}.
$$
Alternatively, one can use the description 
of the Dickson invariants
in terms of Vandermonde-like determinants given in \cite[Prop.~1.3]{wilkerson83}
to verify that $z_n^{n(q-1)-1}$ divides $J$:
$$d_{n,k} = \Delta_1 \Delta_2^{-1} \text{ where }
 \Delta_1 \ = \  
  \det\left(z_j^{q^i}\right)_{\substack{i=0,\ldots n,\ i\neq k\\ j=1,\ldots, n 
    \hphantom{\ j\neq k}}},\quad
 \Delta_2 = {\det}\left(z_j^{q^i}\right)_{\substack{i=0,\ldots n-1\\ j=1,\ldots, n
        \hphantom{-1}}}.
$$
Apply the quotient rule to $\del/\del x_i(d_{n,k})$
and expand $\del/\del x_i(\Delta_1)$ and  $\del/\del x_i(\Delta_2)$
about the $i$-th row.  If $i\neq n$, then $z_n^{q+1}$ divides
 $\del/\del x_i(\Delta_1)$ and  $\del/\del x_i(\Delta_2)$ 
and hence $z_n^{q-1}$ divides $\del/\del x_i (d_{n,k})$.
If $i=n$, expand $\Delta_1$ and $\Delta_2$ about the $n$-th row
and cancel terms to see that 
$z_n^{q-2}$ divides  $\del/\del z_i (d_{n,j})$.  
The last column of this Jacobian matrix is divisible by $z_n^{q-2}$
and the other columns are each divisible by $z_n^{q-1}$.

The Jacobian of the Dickson invariants was also examined by K.~Kuhnigk in her
Doktorarbeit (\cite{kathrin}).
\end{example}

\begin{example} Let $G=\Sl_n(\F_q)$.
As in the last example, let $\F=\F_q$
and $ Q = \prod_{H \in \A} l_H $,
where $\A$ is the set of all hyperplanes in 
$V=\F^n$.
Every reflection in $G$ is a transvection and 
$$
 \F[V]^G \ = \ \F[f_1, f_2, \ldots, f_n],
$$
where $f_1 = Q$ and $f_{i+1}$  is the Dickson invariant $d_{n,i}$
(for $i\geq 1$).
Fix some hyperplane $H$ in $V$ 
and let $G_H$ be its pointwise stabilizer.
In an appropriate basis (with $z_n = l_H$),
$\F[V]^{G_H} = F[u_1, ..., u_n]$
where 
$
   u_i= z_i^q -  z_n^{q-1} z_i$
for $i<n$
and $u_n=z_n$.
Note that the sum of the exponents
of $G_H$ is $\expo = (n-1)(q-1)$.
Since each $f_i$ lies in $\F[V]^{G_H}$,
$f_i = h_i(u_1, \ldots, u_n)$
for some $h_i \in \F[V]$. 
Then 
$   \del f_i/\del z_k = \sum_{j} (\del h_i/\del z_j)(\del u_j /\del z_k)$
is divisible by $z_n^{q-1}$ if $k\neq n$.
Hence, $J=\det\left(\del f_i/\del z_k\right)$ is divisible by $z_n^{(n-1)(q-1)}$ 
and thus
$Q^{\highexpo} = \prod_{H \in \A} l_H^{\highexpo}$ divides $J$.
But $\deg J = (q^n-1)(n-1) =\deg Q(n-1)(q-1)$,
and so 
$$
 J \ =\ Q^{(n-1)(q-1)}\ =\ \prod_{H\in\A}l_H^{\highexpo}.
$$
\end{example}

Note that for both groups $\Gl_n(\F_q)$ and $\Sl_n(\F_q)$,
the Jacobian 
determinant $J$ 
lies in the Hilbert ideal
generated by the basic invariants: 
The image of $J$ in the coinvariant algebra
$\F[V]/(f_1,\ldots, f_n)$ is zero in both cases.

\vspace{2ex}

{\it Acknowledgments}.
We owe thanks to Peter M\"uller and Larry Smith for valuable comments on early
versions of this article, and to Victor Reiner for helpful discussions.
 


\end{document}